\documentclass[12pt]{amsart}
\usepackage{amsmath}
\usepackage{amssymb}
\usepackage{amsmath,amssymb,amsthm,amscd}
\numberwithin{equation}{section}

\def\R{{\mathbb R}}

\newtheorem{theo}{Theorem}[section]
\newtheorem{lemm}{Lemma}[section]
\newtheorem{coro}{Corollary}[section]
\newtheorem{rema}{Remark}[section]

\newtheorem{defi}{Definition}[section]

\def\begeq{\begin{equation}}
\def\endeq{\end{equation}}
\def\p{\partial}

\def\lf{\left}
\def\ri{\right}

\def\R{\Bbb R}

\begin{document}

\title{The Brown-York mass of revolution surface in
asymptotically Schwarzschild manifold}
\author{Xu-Qian Fan$^\dagger $ \&   Kwok-Kun Kwong}
\thanks{$^\dagger $Research partially supported by the National
Natural Science Foundation of China (10901072) and GDNSF
(9451503101004122).}
\address{Department of
Mathematics, Jinan University, Guangzhou,
 510632,
P. R. China.} \email{txqfan@jnu.edu.cn}
\address{Department of
 Mathematics, The Chinese University of Hong Kong,
Shatin, Hong Kong, P. R. China.} \email{kkkwong@math.cuhk.edu.hk}

\renewcommand{\subjclassname}{%
  \textup{2000} Mathematics Subject Classification}
\date{Oct. 2009}
\begin{abstract}
In this paper, we will show that the limit of the Brown-York mass of
a family of convex revolution surfaces in an asymptotically
 Schwarzschild manifold is the ADM mass.
\end{abstract}
\maketitle\markboth{Xu-Qian Fan $\&$   Kwok-Kun Kwong }{Brown-York
mass of revolution surface}

\section{Introduction}
Throughout this paper, we will denote $(\mathbb{R}^3,\delta_{ij})$
as the 3-dimensional Euclidean space, $x^1,x^2,x^3$ as the standard
coordinates, $r$ and $\p$ as the Euclidean distance and the standard
derivative operator on $\R^3$ respectively. Let us first recall some
definitions. First of all, we will adopt the following definition of
asymptotically flat manifolds.

\begin{defi}\label{defaf}
A complete three dimensional manifold $(M,\lambda)$ is said to be asymptotically
flat (AF) of order $\tau$ (with one end) if there is a compact
subset $K$ such that $M\setminus K$ is diffeomorphic to
$\R^3\setminus B_R(0)$ for some $R>0$ and in the standard
coordinates in $\R^3$, the metric $\lambda$ satisfies:
\begin{equation} \label{daf1}
\lambda_{ij}=\delta_{ij}+\sigma_{ij}
\end{equation}
with
\begin{equation} \label{daf2}
|\sigma_{ij}|+r|\p \sigma_{ij}|+r^2|\p\p\sigma_{ij}|+r^3|\p\p\p
\sigma_{ij}|=O(r^{-\tau}),
\end{equation}
for some constant $1\ge\tau>\frac{1}{2}$.
\end{defi}
A coordinate system of $M$ near infinity so that the metric tensor
in this system satisfy the above decay conditions is said to be
admissible. In such a coordinate system, we can define the ADM mass
as follows.
\begin{defi}
The Arnowitt-Deser-Misner (ADM) mass (see \cite{ADM}) of an
asymptotically flat manifold $(M,\lambda)$ is defined as:
\begin{equation} \label{defadm1}
m_{ADM}(M,\lambda)=\lim_{r\to\infty}\frac{1}{16\pi}\int_{S_r}
\lf(\lambda_{ij,i}-\lambda_{ii,j}\ri)\nu^jd\Sigma_r^0,
\end{equation}
where $S_r$ is the Euclidean sphere, $d\Sigma_r^0$ is the volume
element induced by the Euclidean metric, $\nu$ is the outward unit
normal of $S_r$ in $\R^3$ and the derivative is the ordinary partial
derivative.
\end{defi}
We always assume that the scalar curvature is in $L^1(M)$ so that
the limit exists in the definition. In \cite{BTK86}, Bartnik showed that the
ADM mass is a geometric invariant. More precisely,
he proved the following theorem (see \cite[Proposition 4.1]{BTK86} for
a more general setting):
\begin{theo}\label{theb1}
Suppose $(M,\lambda)$ is an  AF manifold with scalar curvature $R(\lambda)\in L^1(M)$.
Let
$\{D_k\}_1^\infty$ be an exhaustion of $M$ by closed sets such that
the set $S_k=\partial D_k$ are connected $C^1$ surfaces without
boundary in $\mathbb{R}^3$ such that
\begin{equation*}
\begin{split}
&r_k=\inf\{|x|,x\in S_k\}\to\infty \text{ as } k\to\infty\\
&r_k ^{-2}\text{Area}(S_k) \text{ is bounded as } k\to\infty.
\end{split}
\end{equation*}
Then
$$m_{ADM}(M,\lambda)=\lim_{k\to\infty}\frac{1}{16\pi}\int_{S_k}
\lf(\lambda_{ij,i}-\lambda_{ii,j}\ri)\nu^jd\Sigma_r^0.$$ That is,
the ADM mass is independent of the sequence of $\{S_k\}$.
\end{theo}
 On the other hand, there have been many studies on the relation between the
 ADM mass of an AF manifold and the so called quasi-local mass. Let us recall the
 definition of the Brown-York quasi-local mass.
 Suppose
$\lf(\Omega, \mu\ri)$ is a compact three dimensional manifold with
smooth boundary $\p \Omega$, if moreover $\p \Omega$ has positive Gauss
curvature, then the Brown-York mass of $\p\Omega$ is
defined as (see \cite{BY1,BY2}):
\begin{defi}
 \begin{equation} \label{defbym1}
m_{BY}\lf(\p \Omega\ri)=\frac{1}{8\pi}\int_{\p \Omega}(H_0-H)d\sigma
\end{equation}
 where  $H$ is the mean curvature of $\p\Omega$ with respect to the
 outward unit normal and
  the metric $\mu$, $d\sigma$ is the volume element induced on $\p\Omega$ by
 $\mu$ and $H_0$ is the mean curvature  of $\p \Omega$ when
  embedded in  $\R^3$.
\end{defi}

The existence of an isometric embedding in $\R^3$ for $\p\Omega$ was
proved by Nirenberg \cite{Niren}, the uniqueness of  the embedding
was given by \cite{Herglotz43,Sacksteder62,PAV}, so the Brown-York
mass is well-defined.

It can be proved that the Brown-York mass and the Hawking quasi-local mass
\cite{SWH} of
the coordinate spheres tends to the ADM mass in some AF manifolds,
see \cite{BY2,HKGH1,BBWYY,BLP,ST02,FST}, and even of nearly round
surfaces \cite{SWW}. It is therefore natural to ask whether the
quasi-local mass of a more general class of surfaces tends to the ADM
mass.

In this paper, we will consider a special class of AF manifolds, called
asymptotically
 Schwarzschild manifolds, which is defined as follows:
\begin{defi}\label{def: metric}
$(N,\tilde{g})$ is called an asymptotically
 Schwarzschild manifold if
   $N=\mathbb{R}^3\setminus K$,  $K$ is a compact set containing the origin, and
  \begin{equation*}
 \tilde{g}_{ij}=\phi^4 \delta_{ij}+b_{ij},
  \phi=\lf(1+\frac{2m}{r}\ri),
    m>0,
  \end{equation*}
  where $|b_{ij}|+r|\p b_{ij}|+r^2|\p\p b_{ij}|
  +r^3|\p\p\p b_{ij}|=O\lf(r^{-2}\ri).$
 \end{defi}
  Clearly, $(N,\tilde{g})$ is an AF manifold. For $b=0$,
  $(N,\tilde{g})$ is called a Schwarzschild
manifold. In this case, we always denote $\tilde{g}$ as $g$. Note
that the scalar curvature of $(N,g)$ is zero \cite{HY} (page 283)
and that of $\lf(N,\tilde{g}\ri)$ is in $L^1(N)$, so in both cases
the ADM mass is well defined.

We will study the limiting behaviors of Brown-York mass on some
family of convex revolution surfaces in an asymptotically
Schwarzschild manifold. We will first prove the following:

\begin{theo}\label{thm1}
Let $(N,{g})$ be a Schwarzschild manifold, $\{D_a\}_{a>0}$ be an
exhaustion of $N$. Suppose $\{S_a=\p D_a\}_{a>0}$ is a family of
closed convex surfaces of revolution in $(\R^3, \delta_{ij})$ with
the rotation axis passing through the origin, satisfying the
following conditions:\\
(i)
\begin{equation}\label{eq: bar K}
\bar K \geq\frac{C_1}{a^2}
\end{equation}
where $\bar K$ is the Gaussian curvature of $S_a$ with induced
Euclidean metric, $C_1>0$ is independent of $a$.\\
(ii)
\begin{equation}\label{eq: bar H}
 0<\bar H\leq  \frac{C_2}{a}
\end{equation}
where $\bar H$ is the mean curvature of $S_a$ with induced
Euclidean metric, $C_2>0$ is independent of $a$.\\
(iii)
\begin{equation}\label{eq:dse}
C_3 a\leq \min_{x\in S_a} r(x)\leq \max_{x\in S_a}r(x)\leq C_4 a
\end{equation}
where $C_i>0$ are independent of $a$. Then
$$\lim_{a\to
\infty}m_{BY}(S_a)=m_{ADM}(N, g).$$
\end{theo}

One example of surfaces satisfying the
conditions in Theorem \ref{thm1} (and also Theorem \ref{thm2} below) is the
family of ellipsoids:
   $$S_a=\lf\{(x^1)^2+(x^2)^2+\frac{(x^3)^2}{4}=a^2\ri\},$$
   which is not nearly round \cite{SWW}. In contrast, the Hawking mass of this family
   of ellipsoids in $(N, g)$ does not tend to the ADM mass of $(N, g)$, indeed one can
   check that the Hawking mass of this family tends to negative infinity as
   $a \rightarrow \infty$.

\begin{rema}\label{rema: prin curv}
The above conditions (i) and (ii) imply that the principal curvature
$\lambda$ of $S_a$ in $(\mathbb{R}^3, \delta)$ satisfy
\begin{equation}
\frac{C_1}{C_2 a}\leq \lambda \leq \frac{C_2}{a}.
\end{equation}
for any $a$. For if $0<\lambda_1\leq \lambda_2$ are the principal
curvatures, then \eqref{eq: bar H} implies
\begin{equation}
\lambda_2\leq \frac{C_2}{a}.
\end{equation}
Together with \eqref{eq: bar K},
\begin{equation}
\lambda_1\geq \frac{C_1}{\lambda_2 a^2}\geq \frac{C_1}{C_2 a}.
\end{equation}
\end{rema}
\begin{rema}
By condition (i) of Theorem \ref{thm1} and the Gauss-Bonnet theorem,
the Euclidean area of $S_a$ is of order $O(a^2)$.
\end{rema}

Then we will generalize the result of Theorem \ref{thm1} to an
asymptotically Schwarzschild manifold:
 \begin{theo}\label{thm2}
 Let $\lf(N,\tilde{g}\ri)$  be an asymptotically Schwarzschild
 manifold and $S$ be a $C^{6,\alpha}$ ($0<\alpha<1$) closed convex revolution surface
 parametrized
 by
 \begeq \label{eq:ass1}
(\bar{w}(\varphi)\cos\theta, \bar{w}(\varphi)\sin\theta,
\bar{h}(\varphi)),\quad 0\leq \theta\leq 2\pi \text{ and } 0\leq
\varphi\leq l.
\endeq
Then there exists $\varepsilon>0$ such that for any family of
$C^{5,\alpha}$ closed convex revolution surfaces $S_a$ in
$\lf(\mathbb{R}^3,\delta\ri)$ satisfying \eqref{eq: bar
K}-\eqref{eq:dse} and is parametrized by
$$\lf(aw_a(\varphi)\cos\theta,\
aw_a(\varphi)\sin\theta,\ ah_a(\varphi)\ri),\  0\leq \theta\leq 2\pi
\text{ and } 0\leq \varphi\leq l$$ such that \begeq \label{eq:ass2}
|w_a-\bar{w}|_{C^{4}}+|h_a-\bar{h}|_{C^{4}}\leq \varepsilon\quad
\text{ for } a>>1,
\endeq
we have
$$\lim_{a\to\infty}m_{BY}(S_a)=m_{ADM}(N, \tilde g).$$

\end{theo}
From this result, one has
 \begin{coro} \label{co1}
 Let $\lf(N,\tilde{g}\ri)$  be an asymptotically Schwarzschild
 manifold.
Let $\{S_i\}$ be a family of $C^7$ closed convex revolution surfaces
 in $\lf(\mathbb{R}^3,\delta\ri)$ satisfying \eqref{eq: bar K}-\eqref{eq:dse}
 and is
parametrized as:
$$\lf(a_iw_i(\varphi)\cos\theta,\
a_iw_i(\varphi)\sin\theta,\ a_ih_i(\varphi)\ri),\  0\leq \theta\leq
2\pi \text{ and } 0\leq \varphi\leq l$$ for some constant $l>0$,
 here $a_i$ are positive numbers with $\displaystyle\lim_{i\to\infty}a_i=+\infty.$
 If
 there is a constant $c$ such that
  \begeq \label{eq:ass3}
\ |w_i|_{C^{7}}+|h_i|_{C^{7}}\leq c
\endeq for all $i$,
then there is a subsequence $\{S_{i_k}\}$ of $\{S_i\}$ such that
$$\lim_{k\to\infty}m_{BY}\lf(S_{i_k}\ri)=m_{ADM}(N, \tilde g).$$
 \end{coro}

This paper is organized as follows. In section 2, we will prove
Theorem \ref{thm1}. Theorem \ref{thm2} will be proved in section 3.

\section*{Acknowledgments}
The authors would like to thank Prof. Luen-fai Tam for his inspired
guidance, constant encouragement and very worthy advice. The author
 (XQ) would also like to thank Prof. Youde Wang for his invitation.

\section{Proof of Theorem \ref{thm1}}
For simplicity, from now on to the end of the paper, we use
$O\lf(a^{-1}\ri)$ to denote a quantity which is
bounded by $Ca^{-1}$ for some constant $C$ independent of $a$ as $a$
is big enough. Similar for $O\lf(1\ri)$, $O\lf(a^{-2}\ri)$ and
$O\lf(a^{-3}\ri)$. We will first compute the mean curvature of $S_a$
in $(N,{g})$ and of the embedded surface of the Euclidean space
respectively.

From the assumptions of $S_a$, we can assume that $S_a$ is
parametrized by
\begin{equation}
(a w_a(\varphi)\cos \theta, aw_a(\varphi)\sin \theta,
ah_a(\varphi)), \quad 0\leq \varphi\leq l_a, 0\leq \theta\leq
2\pi,
\end{equation}
$w_a(\varphi), h_a(\varphi)$ being smooth functions for $\varphi\in
[0, l_a]$ (i.e. $w_a$, $h_a$ can be extended smoothly on
a slightly larger interval ). Moreover,\\
(i)\begin{equation}\label{eq: wh}
\begin{split}
h_a(0)>h_a(l_a)\\
C_3\leq \sqrt{w_a^2+h_a^2}\leq C_4 \\
w_a>0 \text{ on $(0, l_a)$},
\end{split}
\end{equation}
\noindent (ii) The generating curve $(w_a(\varphi), h_a(\varphi))$
is parameterized by arc length. i.e.
\begin{equation}\label{eq: arclength}
w_a'^2+h_a'^2=1.
\end{equation}

\noindent (iii) $w_a$ is anti-symmetric about $0$ and $l_a$, $h_a$
is symmetric about $0$ and
$l_a$, i.e. \\
\begin{equation}\label{eq: odd}
\begin{split}
w_a(-\varphi)&=-w_a(\varphi),  \quad w_a(l_a+\varphi)=-w_a(l_a-\varphi),\\
h_a(-\varphi)&=h_a(\varphi),  \quad h_a(l_a+\varphi)=h_a(l_a-\varphi) \text{ for
$\varphi\in [0,\varepsilon) $.}
\end{split}
\end{equation}
This implies
\begin{equation}\label{eq: w0}
w_a(0)=w_a(l_a)=h_a'(0)=h_a'(l_a)=0.
\end{equation}
Since $S_a$ is convex in $(\mathbb{R}^3,\delta)$ and the Gaussian
curvature $\bar K$ of $S_a$ with the induced metric $d\bar{s}^2$ is
\begin{equation}
\bar K=\frac{h_a'(w_a'h_a''-w_a''h_a')}{a^2w_a}\quad \text{for $\varphi\in
(0,l_a)$.}
\end{equation}
So $h_a'<0$ for $\varphi\in (0, l_a)$ by \eqref{eq: wh}.

Let $\phi_a$ be the function $\phi$ restricted on $S_a$, note that
in $(\varphi, \theta)$ coordinates, $\phi_a=\phi_a(\varphi)$ is
independent of $\theta$. We have the following lemma:

\begin{lemm}\label{lem: extend}
The functions $\displaystyle \frac{w_a}{h_a'}$ and $\displaystyle
\frac{\phi_a'}{h_a'}$ can be extended continuously to the whole $[0,
l_a]$. Moreover there exists a constant $C$ independent of $a$ such
that
\begin{equation}
\lf|\frac{w_a}{h_a'}\ri|\leq C, \quad \lf|\frac{\phi_a'}{h_a'}\ri|\leq \frac{C}{a}
\end{equation}
for all $a$.
\end{lemm}

\begin{proof}
We first show that the limits
\begin{equation}
\lim_{\varphi\to 0}\frac{w_a}{h_a'}, \lim_{\varphi\to l_a}\frac{w_a}{h_a'},
\end{equation}
exist and are bounded by $C$.

The Gaussian curvature $\bar K$ of the point $(0, 0, ah_a(0))$ on
$S_a$ with induced Euclidean metric is equal to
\begin{equation}
\bar K=\frac{h_a''(0)^2}{a^2}.
\end{equation}
This can be shown by the fact that for an arc-length parametrized
plane curve $(w_a(\varphi), h_a(\varphi))$, its curvature is given
by $-w_a''h_a'+h_a''w_a'$. So at $(0, h_a(0))$, its curvature is
$h_a''(0)$.

As $\displaystyle\bar K\geq\frac{C_1}{a^2}$ by \eqref{eq: bar K},
$|h_a''(0)|\geq \sqrt{C_1}>0$. By L'Hospital rule,
\begin{equation}
\lim_{\varphi\to 0}\frac{w_a}{h_a'}=\frac{w_a'(0)}{h_a''(0)}
\end{equation}
which is finite and is bounded by some $C>0$ by \eqref{eq: bar K}
and \eqref{eq: arclength}. The same applies to $\displaystyle
\lim_{\varphi\to l_a}\frac{w_a}{h_a'}$.

Next, observe that one of the principal curvatures of $S_a$ in
$(\mathbb{R}^3, \delta)$ is $\displaystyle -\frac{h_a'}{aw_a}$ (
\cite{doCarmo} p.162, (10)). So by Remark \ref{rema: prin curv}, we
have
\begin{equation}
\lf| \frac{w_a}{h_a'}\ri|\leq C
\end{equation}
on the whole $[0, l_a]$ for all $a$.

By differentiating $\displaystyle
\phi_a=1+\frac{m}{2a\sqrt{w_a^2+h_a^2}}$,
\begin{equation}
\frac{\phi_a'}{h_a'}=-\frac{m}{2a(w_a^2+h_a^2)^{\frac{3}{2}}}
(w_a'\frac{w_a}{h_a'}+h_a)
\end{equation}
which can be extended to $[0, l_a]$ by the above, and is of order
$O\lf(a^{-1}\ri)$ by \eqref{eq: wh}, \eqref{eq: arclength}.
\end{proof}
We have the following estimates
\begin{lemm}\label{lem: phi'}
Regarding $\phi_a=\phi_a(\varphi)$ as functions on $S_a$, we have
\begin{equation}
\phi_a'=O(a^{-1}), \quad \phi_a''=O(a^{-1}).
\end{equation}
\end{lemm}

\begin{proof}
Let $A=w_a^2+h_a^2$. As $\displaystyle \phi_a=1+\frac{m}{2a
\sqrt{A}}$, we only have to prove
\begin{equation}
(A^{-\frac{1}{2}})'=O(1), \quad (A^{-\frac{1}{2}})''=O(1).
\end{equation}
By direct computations and \eqref{eq: wh}, \eqref{eq: arclength},
\begin{equation}
\begin{split}
|(A^{-\frac{1}{2}})'|&=|A^{-\frac{3}{2}}(w_a w_a'+h_a h_a')|\\
&\leq A^{-\frac{3}{2}}(w_a^2+h_a^2)^{\frac{1}{2}}(w_a'^2+h_a'^2)^{\frac{1}{2}}\\
&=O(1)
\end{split}
\end{equation}
and
\begin{equation}
\begin{split}
|(A^{-\frac{1}{2}})''|
&=\lf|\frac{3}{2}A^{-\frac{5}{2}}(w_aw_a'+h_a h_a')^2-A^{-\frac{3}{2}}
(1+w_a w_a''+h_a h_a'')\ri|\\
&\leq \frac{3}{2}A^{-\frac{5}{2}}(w_a^2+h_a^2)+A^{-\frac{3}{2}}
(1+(w_a^2+h_a^2)^{\frac{1}{2}} (w_a''^2+h_a''^2)^{\frac{1}{2}}).
\end{split}
\end{equation}
The two principal curvatures of $S_a$ with induced Euclidean metric
are $\displaystyle-\frac{h_a'}{aw_a}$ and
$a^{-1}(w_a''^2+h_a''^2)^{\frac{1}{2}}$ (\cite{doCarmo} p.162,
(10)), hence by Remark \ref{rema: prin curv},
$|(A^{-\frac{1}{2}})''|=O(1). $
\end{proof}

 From now on, we will drop the
subscript $a$ and denote $w_a$ by $w$, $h_a$ by $h$, $\phi_a$ by
$\phi$ and $l_a$ by $l$. We also denote $ds^2$ to be the metric on
$S_a$ induced from $g$.

\begin{lemm}\label{lem:pgc}
 The Gaussian curvature ${K}$ of $(S_a, ds^2)$ is positive for large enough $a$.
 In particular, there exists a unique isometric
 embedding of $(S_a, ds^2)$ into $(\mathbb{R}^3, \delta)$ for
 sufficiently large $a$.
\end{lemm}
\begin{proof}
Let $d\bar s^2$ be the metric on $S_a$ induced by $\delta$ and $d
s^2$ be that induced by $g$. Since
\begin{equation*}
\begin{split}
d\bar{s}^2&=a^2(d\varphi^2+w^2d\theta^2)=\bar{E}d\varphi^2+\bar{G}d\theta^2
\text{ and }\\
d {s}^2&= {\phi}^4d\bar{s}^2= {E}d\varphi^2+ {G}d\theta^2,
\end{split}
\end{equation*}
one has
\begin{equation*}
\begin{split}
 {E}=\bar{E}+O(a), & {E}_\varphi=\bar{E}_\varphi+O(a),
 {E}_{\varphi\varphi}=\bar{E}_{\varphi\varphi}+O(a)\\
\text{ and   } \qquad & {E}_\theta=\bar{E}_\theta+O(a),
 {E}_{\theta\theta}=\bar{E}_{\theta\theta}+O(a).
\end{split}
\end{equation*}
Similar result holds for $ {G}.$ By the formula
\begin{equation}
 K=-\frac{1}{\sqrt{ E  G}}\lf(\lf(\frac{ E_\theta}{\sqrt{ E  G}}\ri)_\theta
 + \lf(\frac{ G_\phi}{\sqrt{ E  G}}\ri)_\phi\ri)
\end{equation}
and the corresponding formula for $\bar  K$, one can get $ {K}=\bar
K+O\lf(a^{-3}\ri)$. Hence the lemma holds.
\end{proof}

 Now let us compute the mean curvature of a revolution surface in
$(\mathbb{R}^3,\delta)$.
\begin{lemm}\label{lem: bar H}
For a smooth revolution surface $S$ in $(\mathbb{R}^3,\delta)$
parametrized by
\begin{equation}
(au(\varphi)\cos\theta, au(\varphi)\sin\theta, av(\varphi)), \quad
0<\varphi<l, 0<\theta< 2\pi,
\end{equation}
its mean curvature $\bar{H}$ with respect to $\delta$ is
\begin{equation}
\bar
H=\frac{u''}{aTv'}-\frac{T'u'}{aT^2v'}-\frac{v'}{aTu}\quad\text{where
$T=\sqrt{u'^2+v'^2}$.}
\end{equation}
\end{lemm}

\begin{proof}
The mean curvature $\bar{H}$ of $S$ with respect to $\delta$ is
computed to be
\begin{equation}
\bar{H}=\frac{v'u''-u'v''}{aT^3}-\frac{v'}{a Tu}.
\end{equation}
By differentiating $T^2$,
\begin{equation}
u'u''+v'v''=TT'
\end{equation}
which implies
\begin{equation}
\begin{split}
v'u''-u'v''&= v'u''+\frac{u'^2u''-u'TT'}{v'}\\
&=\frac{(u'^2+v'^2)u''-u'TT'}{v'}\\
&=\frac{T^2 u''-TT'u'}{v'}.
\end{split}
\end{equation}
So we have
\begin{equation*}
\bar H=\frac{u''}{aTv'}-\frac{T'u'}{aT^2v'}-\frac{v'}{aTu}.
\end{equation*}
The lemma is proved.
\end{proof}

\begin{lemm}\label{lem: H}
The mean curvature ${H}$ of $S_a$ with respect to ${g}$ is
\begin{equation}\label{eq:tH}
{H}=\frac{w''}{a\phi^2h'}-\frac{h'}{a\phi^2 w}+4\phi^{-3}n(\phi)
\end{equation}
where $n$ is the outward unit normal vector of $S_a$ with respect to
$\delta$.
\end{lemm}

\begin{proof}
By Lemma \ref{lem: bar H}, the mean curvature of $S_a$ with respect
to $\delta$ is
\begin{equation}
\bar H=\frac{w''}{ah'}-\frac{h'}{a w}.
\end{equation}
The mean curvature $ {H}$ of $S_a$ with respect to $ {g}$ is
(\cite{SYL79}, page 72):
\begin{equation}
 {H}= {\phi}^{-2}\lf(\bar{H}+4 {\phi}^{-1}n\lf( {\phi}\ri)\ri)
\end{equation}
where $n$ is the outward unit normal vector of $S_a$ with respect to $\delta$. \\
Therefore
\begin{equation*}
 {H}=\frac{w''}{a {\phi}^2h'}-\frac{h'}{a {\phi}^2w}+4 {\phi}^{-3}
n\lf( {\phi}\ri).
\end{equation*}
\end{proof}

\begin{lemm}\label{lem: embed}
For sufficiently large $a$, there is an isometric embedding of
$(S_a,ds^2)$ into $(\mathbb{R}^3,\delta)$ which is given by
\begin{equation}\label{eq: embed}
x^1=au(\varphi)\cos\theta,\ x^2=au(\varphi)\sin\theta,\
x^3=av(\varphi), \ 0\leq \varphi\leq l, 0\leq \theta\leq 2\pi
\end{equation}
where
\begin{equation}
\begin{split}
u&=\phi^2w\\
v'&=\phi^2h' \lf(1-\frac{2\phi'ww'}{h'^2}+O\lf(a^{-2}\ri)\ri)\\
u'^2+v'^2&=\phi^4.
\end{split}
\end{equation}
\end{lemm}

\begin{proof}
The existence has already been proved in Lemma \ref{lem:pgc}.

In $(\varphi, \theta)$ coordinates, the metric on $S_a$ induced by
$g$ can be written as:
\begin{equation}\label{eq: metric}
ds^2=a^2\phi^4d\varphi^2+a^2\phi^4 w^2 d\theta^2.
\end{equation}
We can regard $(S_a,ds^2)$  as $\mathbb{S}^2$, the sphere with the
metric $ds^2$. Now we want to find two functions $u,v$ such that the
surface written as  the form \eqref{eq: embed} is an embedded
surface  $S^e_a$ of $S_a$ into $(\mathbb{R}^3,\delta)$. First of
all, the induced metric by the Euclidean metric on the surface which
is of the form \eqref{eq: embed}  can be written as:
$$ds^2_e=a^2\lf(u'^2+v'^2\ri)d\varphi^2+a^2u^2d\theta^2 .$$
Comparing this with \eqref{eq: metric}, one can choose
\begin{equation}\label{eq: u=phi^2w}
u=\phi^2w .
\end{equation}
Consider
\begin{equation}\label{eq: phi^4-u'^2}
\begin{split}
\phi^4-u'^2
&=\phi^2 (\phi^2-(2\phi'w+\phi w')^2)\\
&=\phi^2 (\phi^2 (w'^2+h'^2)-(2\phi'w+\phi w')^2)\\
&=\phi^2 (\phi^2h'^2-4\phi\phi'ww'-4\phi'^2w^2)\\
&=\phi^4h'^2 \lf(1-\frac{4\phi'ww'}{\phi
h'^2}-\frac{4\phi'^2w^2}{\phi^2 h'^2}\ri) .
\end{split}
\end{equation}
By Lemma \ref{lem: extend} and Lemma \ref{lem: phi'}, the functions
$\displaystyle\frac{\phi'ww'}{\phi h'^2}, \frac{\phi'^2w^2}{\phi^2
h'^2}$ can be extended continuously on $[0, l]$ with
$\displaystyle\frac{\phi'ww'}{\phi h'^2}=O(a^{-1}),
\frac{\phi'^2w^2}{\phi^2 h'^2}=O(a^{-2})$. So $\displaystyle
1-\frac{4\phi'ww'}{\phi h'^2}-\frac{4\phi'^2w^2}{\phi^2 h'^2}> 0$
for sufficiently large $a$. For these $a$, we can take
\begin{equation}
v'=\phi^2 h' \lf(1-\frac{4\phi'ww'}{\phi
h'^2}-\frac{4\phi'^2w^2}{\phi^2 h'^2}\ri)^{\frac{1}{2}} ,
\end{equation}
so that
\begin{equation}
u'^2+v'^2=\phi^4.
\end{equation}
Note that by \eqref{eq: odd}, $v'$ is an odd function for
$\varphi\in[-l,l].$ By choosing an initial value, one can get an
even function $v$. By the above argument, one has
\begin{equation}
v'=\phi^2h' \lf(1-\frac{2\phi'ww'}{ h'^2}+O\lf(a^{-2}\ri)\ri).
\end{equation}

From \eqref{eq: u=phi^2w} and \eqref{eq: phi^4-u'^2}, near
$\varphi=0$, $u,\ v$ can be extended naturally to
$(-\varepsilon,\varepsilon)$ for some $\varepsilon>0$. Since $u$ is
an odd function in $\varphi$ , $v$ is an even function in $\varphi$,
and $u'^2+v'^2=T^2>0$, the generating curve in $\{x^2=0\}$ is
symmetric with respect to $x^3$-axis, and is smooth at $\varphi=0$.
Similarly, it is also smooth at $\varphi=l$. Hence the revolution
surface determined by the choice of $u,\ v$ as above, can be
extended smoothly to a closed revolution surface, which is an
embedded surface of $S_a$. This completes the proof of the lemma.
\end{proof}

 Now we
are ready to prove Theorem \ref{thm1}.
\begin{proof}[Proof of Theorem \ref{thm1}]
Let $u, v$ be defined as in Lemma \ref{lem: embed}. Recall that
\begin{equation}\label{eq: u}
\begin{split}
u&=\phi^2w\\
v'&=\phi^2h' \lf(1-\frac{2\phi'ww'}{h'^2}+O\lf(a^{-2}\ri)\ri)\\
u'^2+v'^2&=\phi^4=T^2 \quad\text{ where $T=\phi^2$.}
\end{split}
\end{equation}
By Lemma \ref{lem: phi'}, we have
\begin{equation}\label{eq: u''}
\begin{split}
T'&=2\phi'+O\lf(a^{-2}\ri)\\
u'&=\phi^2w'+O\lf(a^{-1}\ri)\\
u''&=\phi^2w''+4\phi' w'+2\phi''w+O\lf(a^{-2}\ri).
\end{split}
\end{equation}
By Lemma \ref{lem: bar H} and Lemma \ref{lem: embed},
\begin{equation}\label{eq: H0}
H_0=\frac{u''}{aTv'}-\frac{T'u'}{aT^2v'}-\frac{v'}{aT u}.
\end{equation}

 \noindent
 Combining with Lemma \ref{lem: H},
\begin{equation}\label{eq: H0-H}
 {H}_0- {H}=\lf(\frac{u''}{aTv'}-\frac{w''}{a\phi^2
h'}\ri)-\frac{T'u'}{aT^2v'}-\lf(\frac{v'}{aTu}-\frac{h'}{a\phi^2
w}\ri)-4\phi^{-3}n(\phi).
\end{equation}
Using \eqref{eq: u} and \eqref{eq: u''},
\begin{equation}\label{eq: I}
\begin{split}
\frac{u''}{aTv'}-\frac{w''}{a\phi^2 h'}
&=\frac{w''}{a\phi^2  h'}+\frac{4\phi' w'}{a h'}+\frac{2\phi''w}{a h'}+
\frac{2\phi'ww'w''}{a h'^3}-\frac{w''}{a\phi^2 h'}+O\lf(a^{-3}\ri)\\
&=\frac{4\phi' w'}{a h'}+\frac{2\phi''w}{a h'}+\frac{2\phi'ww'w''}{a
h'^3}+O\lf(a^{-3}\ri).
\end{split}
\end{equation}
By \eqref{eq: u} and \eqref{eq: u''},

\begin{equation}\label{eq: II}
\begin{split}
-\frac{T'u'}{aT^2v'} &=-\frac{2\phi'w'}{ah'}+O\lf(a^{-3}\ri).
\end{split}
\end{equation}
By \eqref{eq: u},
\begin{equation}\label{eq: III}
\begin{split}
-\frac{v'}{aTu}+\frac{h'}{a\phi^2 w}
&=-\frac{h'}{a\phi^2 w}+\frac{2\phi'w'}{a h'}+\frac{h'}{a\phi^2 w}+O\lf(a^{-3}\ri)\\
&=\frac{2\phi'w'}{a h'}+O\lf(a^{-3}\ri).
\end{split}
\end{equation}
Summing \eqref{eq: I}, \eqref{eq: II} and \eqref{eq: III} and
comparing with \eqref{eq: H0-H}, we have

\begin{equation}
 {H}_0- {H} =\frac{4\phi' w'}{a h'}+\frac{2\phi''w}{a
h'}+\frac{2\phi'ww'w''}{a h'^3}-4\phi^{-3}n(\phi)+O\lf(a^{-3}\ri).
\end{equation}
As $w'w''=-h'h''$ by \eqref{eq: arclength}, so
\begin{equation}
 {H}_0- {H} =\frac{4\phi' w'}{a h'}+\frac{2\phi''w}{a
h'}-\frac{2\phi'ph''}{a h'^2}-4\phi^{-3}n(\phi)+O\lf(a^{-3}\ri).
\end{equation}
Denote the Euclidean area element of $S_a$ by $d\sigma_0$, the area
element of $(S_a,d {s}^2)$ by $d\sigma$. Note that $ {H}_0-
{H}=O\lf(a^{-2}\ri)$, $d\sigma-d\sigma_0=O\lf(a^{-1}\ri)d\sigma_0$
and $\displaystyle \int_{S_a}d\sigma_0=O\lf(a^2\ri)$. To prove the
result, it suffices to show
$$\lim_{a\to \infty} \frac{1}{8\pi} \int_{S_a}\lf( {H}_0
- {H}\ri)d\sigma_0=m.$$
The Euclidean area element is computed to be
\begin{equation}
d\sigma_0 =a^2 w d\varphi d\theta.
\end{equation}
By \eqref{eq: w0} and Lemma \ref{lem: extend},
\begin{equation}
\begin{split}
&\int_{S_a} (\frac{4\phi' w'}{a h'}+\frac{2\phi''w}{a
h'}-\frac{2\phi' w h''}{a h'^2}) d\sigma_0\\
&=2\pi a\int_{0}^l (\frac{4\phi' ww'}{h'}+\frac{2\phi''w^2}{h'}
-\frac{2\phi'w^2h''}{ h'^2}) d\varphi\\
&=2\pi a\int_{0}^l  \frac{d}{d\varphi}\bigg(\frac{2\phi' w^2}{h'}\bigg)d\varphi\\
&=0.
\end{split}
\end{equation}
Since the norm of the Euclidean gradient of $\phi$ has $|\nabla _0
\phi|=O(r^{-2})$, one has $n(\phi)=O(a^{-2})$. So
\begin{equation}
\begin{split}
\frac{1}{8\pi}\int_{S_a}\lf( {H}_0- {H}\ri)d\sigma_0
&=-\frac{1}{2\pi}\int_{S_a} \phi^{-3}n(\phi)d\sigma_0+O\lf(a^{-1}\ri)\\
&=-\frac{1}{2\pi}\int_{S_a} n(\phi)d\sigma_0+O\lf(a^{-1}\ri).
\end{split}
\end{equation}
By Theorem \ref{theb1}, or Proposition 4.1 in \cite{BTK86}, the
definition of the ADM mass of $N$ can be taken as
\begin{equation}\label{eq:admn1}
\lim_{a\to \infty}\frac{1}{16\pi} \int_{S_a} \sum_{i,
j}( {g}_{ij,i}- {g}_{ii,j})n^j d\sigma_0=m.
\end{equation}
where $n$ is the unit outward normal of $S_a$ with respect to
$\delta$. By a direct computation,
\begin{equation}\label{eq:admn2}
\sum_{i, j}( {g}_{ij,i}- {g}_{ii,j})n^j=-8 {\phi}^3
n^j\frac{\p {\phi}}{\p x^j}=-8 n(\phi)+O\lf(a^{-3}\ri).
\end{equation}
Combining \eqref{eq:admn1} and \eqref{eq:admn2}, we have
\begin{equation}
m=-\lim_{a\to \infty}\frac{1}{2\pi}\int_{S_a} n(\phi)d\sigma_0.
\end{equation}
Therefore
\begin{equation}
\lim_{a\to \infty} \frac{1}{8\pi}
\int_{S_a}( {H}_0- {H})d\sigma=\lim_{a\to \infty}
\frac{1}{8\pi} \int_{S_a}( {H}_0- {H})d\sigma_0=m.
\end{equation}
We are done.
\end{proof}

\section{Proof of Theorem \ref{thm2}}
We will reduce the case of Theorem \ref{thm2} to the
 Schwarzschild manifold. Let us compare the mean curvature before
 embedding.
 \begin{lemm}\label{lmchb}
For the surfaces $S_a$ satisfying the conditions in Theorem
\ref{thm1}, we have
$$|\tilde H - H|\leq Ca^{-3}$$
for some constant $C$ independent of $a$, where $\tilde H$ and $H$
are the mean curvatures of $S_a$ with respect to $\tilde g$ and $g$
respectively.
\end{lemm}

\begin{proof}
We claim that
\begin{equation}\label{eq: tilde A-A}
|\tilde A-A|_g=O\lf(a^{-3}\ri)
\end{equation}
where $A$ and $\tilde A$ are the second fundamental forms with
respect to $g$ and $\tilde g$ respectively.

Let $\rho(x)$ defined on $N$ to be the distance from $x$ to $S_a$
with respect to $\tilde g$. We will use the fact \cite[(7.10)]{HI}:
\begin{equation}\label{eq: A}
\tilde A(X, Y)-|\nabla \rho|_g A(X, Y)=\lf(\Gamma_{ij}^k-\tilde
\Gamma_{ij}^k\ri)X^iY^j \rho_k
\end{equation}
for any tangent vectors $X, Y$ of $S_a$. For completeness, we prove
it here. We proceed as in \cite{SWW} Lemma 2.6. First of all, we
have
\begin{equation}\label{eqsedwg}
\begin{split}
A(X, Y)&=g\lf(\nabla_X \left(\frac{\nabla \rho}{|\nabla \rho|_g}\right), Y\ri)\\
&=\frac{g(\nabla_X (\nabla \rho), Y)}{|\nabla \rho|_g}\\
&=\frac{ X(Y(\rho))-(\nabla_X Y)(\rho)}{|\nabla \rho|_g}\\
&=\frac{ X^iY^j\rho_{ij}-X^iY^j\Gamma_{ij}^k \rho_k}{|\nabla
\rho|_g},
\end{split}
\end{equation}
here the subscript denotes ordinary derivative and $\Gamma_{ij}^k$
are the Christoffel symbols with respect to $g$, with the indices
$i, j, k=1, 2, 3$. Denote $\tilde \Gamma_{ij}^k$ to be the
Christoffel symbols with respect to $\tilde g$. Then since the
$\tilde g$ gradient $|\tilde \nabla \rho|_{\tilde g}=1$, we also
have
\begin{equation}
\tilde A(X, Y)= X^iY^j\rho_{ij}-X^iY^j\tilde\Gamma_{ij}^k \rho_k.
\end{equation}
Combining this with \eqref{eqsedwg}, we can get \eqref{eq: A}.

Note that $|\Gamma_{ij}^k-\tilde \Gamma_{ij}^k|=O\lf(r^{-3}\ri)$ by
the assumptions of the metrics. By asymptotic flatness,
$\displaystyle 1=\tilde g^{ij}\rho_i\rho_j\geq C \sum \rho_i^2$, so
$|\rho_i|$ is uniformly bounded. The condition $\tilde
g_{ij}=g_{ij}+b_{ij}$ implies $|\tilde g^{ij}-g
^{ij}|=O\lf(r^{-2}\ri)$, so
\begin{equation}
||\nabla \rho|_g^2-1|=|(g^{ij}-\tilde
g^{ij})\rho_i\rho_j|=O\lf(r^{-2}\ri)
\end{equation}
which implies
\begin{equation}
|\nabla \rho|_g=1+O\lf(r^{-2}\ri).
\end{equation}
Finally, the principal curvatures $\bar \lambda_i$ in Euclidean
metric are of order $O\lf(a^{-1}\ri)$ by Remark \ref{rema: prin
curv}, the principal curvatures $\lambda_i$ with respect to $g$ are
related to $\bar \lambda_i$ by (\cite{HY} Lemma 1.4):
\begin{equation}
\lambda_i=\phi^{-2}\bar \lambda_i+2\phi^{-3}n(\phi)
\end{equation}
where $n$ is the unit outward normal with respect to $\delta$. In
particular, as $n(\phi)=O(a^{-2})$,
\begin{equation}
|A|_g=O(a^{-1}).
\end{equation}
Combining all these together with \eqref{eq: A}, it is easy to see
that \eqref{eq: tilde A-A} holds. Combining \eqref{eq: tilde A-A}
and the metric conditions of $g$ and $\tilde g$ in Definition
\ref{def: metric}, this implies the lemma.
\end{proof}

Let $\lf(S_a, d\tilde{s}^2\ri),\ \lf(S_a, ds^2\ri)$ denote the
surface $S_a$ with metric $d\tilde{s}^2, ds^2$ induced from
$\tilde{g},\ g$ respectively. By Lemma \ref{lem:pgc}, for $a>>1$,
the Gaussian curvatures on $\lf(S_a, d\tilde{s}^2\ri)$ and $\lf(S_a,
ds^2\ri)$ are both positive, which implies that they can be
isometrically embedded into $\lf(\mathbb{R}^3,\delta\ri)$ uniquely.
Now let us compare the mean curvature after embedding:
\begin{lemm}\label{lmcha}
Under the same notations and conditions of Theorem \ref{thm2}. Let
$\tilde{H}_0,\ H_0$ be the mean curvature of the embedded surfaces
of $\lf(S_a, d\tilde{s}^2\ri)$ and $\lf(S_a, ds^2\ri)$ in
$(\mathbb{R}^3, \delta)$ respectively, as $a>>1$, we have
$|\tilde{H}_0-H_0|\leq C_5a^{-3}$ for some constant $C_5(S)$.
\end{lemm}
\begin{proof}
We can set $\hat{\varphi}=\frac{\pi}{l}\varphi$, so it suffices to
show that the lemma holds for $l=\pi$. Also, by identifying $S$ and
$S_a$ with the sphere $\mathbb{S}^2$, we can regard all the metrics
here ($ds^2$ etc.) to be metrics on $\mathbb{S}^2$. We will denote
$w_a$ as $w$ and $h_a$ as $h$. Similar to \eqref{eq: metric}, one
has
 \begin{equation}
 \begin{split}
  d\bar{s}^2&=a^2\lf(\lf(\lf(w'\ri)^2+\lf(h'\ri)^2\ri)d\varphi^2+w^2d\theta^2\ri)
  \text{ and }\\
d\bar{s}^2_S&=\lf(\lf(\bar{w}'\ri)^2+\lf(\bar{h}'\ri)^2\ri)d\varphi^2
+\bar{w}^2d\theta^2
\end{split}
  \end{equation}
which are the metrics on $S_a$ and $S$ induced from the Euclidean
metric respectively. By definition,
 \begin{equation}
 \begin{split}
ds^2&=\phi^4d\bar{s}^2,\
d\tilde{s}^2=ds^2+b_{ij}dx^idx^j\big|_{S_a}.
\end{split}
\end{equation}

From \eqref{eq:ass2}, $w$ and its derivatives up to forth order are
uniformly bounded for $a>>1$, the same holds for $h$. By the
conditions of $b_{ij}$, it is easy to see that the followings hold:
\begin{equation} \label{12cmtg}
\begin{split}
\|a^{-2} d\tilde{s}^2- a^{-2} ds^2\|_{C^3}=a^{-2}\|b_{ij}dx^idx^j|_{S_a}\|_{C^3}
\leq C_6a^{-2},
\end{split}
\end{equation}
\begin{equation}\label{eq: g-bar g}
\|a^{-2} ds^2- a^{-2} d\bar{s}^2\|_{C^3}=a^{-2}\|(\phi^4-1)d\bar{s}^2\|_{C^3}
\leq C_6a^{-1}
\end{equation}
for some constant $C_6(S)$. By \eqref{eq:ass2}, we have
 \begin{equation} \label{12cmeg}
 \|a^{-2}d\bar{s}^2-d\bar{s}^2_S\|_{C^3}\leq
C_7\varepsilon
 \end{equation}
 for some constant $C_7(S)$.  So for $a>>1$, by \eqref{eq: g-bar g} and
 \eqref{12cmeg}, we have
 $$\|a^{-2}ds^2-d\bar{s}^2_S\|_{C^3}\leq
\lf(C_6+C_7\ri)\varepsilon.$$
 By the result of \cite[Lemma 5.3]{MST}, if we choose some
 $\displaystyle 0<\varepsilon<\frac{\delta}{\pi^{1-\alpha}(C_6+C_7)}$
 such that
$$\|
a^{-2}ds^2-d\bar{s}^2_S\|_{C^{2,\alpha}}<\delta$$ for $a$ big
enough, where $\delta$ is the one given by \cite[Lemma 5.3]{MST},
then there are isometric embeddings $\tilde{X}$ and $X$ of
$\lf(\mathbb{S}^2,a^{-2}d\tilde{s}^2\ri)$ and
$\lf(\mathbb{S}^2,a^{-2}ds^2\ri)$ respectively, such that by
\eqref{12cmtg}, for sufficiently large $a$,
$$\|\tilde{X}-X\|_{C^{2,\alpha}}\leq C_8\|a^{-2} d\tilde{s}^2
- a^{-2} ds^2\|_{C^{2,\alpha}}=O\lf(a^{-2}\ri)$$
for some constant $C_8(S)$. Since $a\tilde{X}, aX$ are the isometric
embeddings of $\lf(\mathbb{S}^2,d\tilde{s}^2\ri)$ and
$(\mathbb{S}^2,ds^2)$ respectively. Hence
$|\tilde{H}_0-H_0|=O\lf(a^{-3}\ri).$ The lemma holds.
\end{proof}

Now we can prove Theorem \ref{thm2}.
\begin{proof}[Proof of Theorem \ref{thm2}]
By Theorem \ref{thm1}, we know that
$$\lim_{a\to\infty}\frac{1}{8\pi}\int_{S_a}\lf(H_0-H\ri)d\sigma=m_{ADM}(N, g).$$
Since the ADM mass of $(N,g)$ is equal to that of
$\lf(N,\tilde{g}\ri)$, combining with Lemma \ref{lmchb} and Lemma
\ref{lmcha}, we can get the result.
\end{proof}

\end{document}